\documentclass[10pt]{amsart}
\usepackage{amsmath,amsfonts,amssymb,amscd,bm,hyperref}

\newcommand{\al}{\alpha}    \newcommand{\be}{\beta}
    
\newcommand{\ep}{\epsilon}  
\newcommand{\la}{\lambda}

\def\<{\langle}             \def\>{\rangle}
\newcommand{\R}{\mathbb{R}}\newcommand{\Z}{\mathbb{Z}}
\newcommand{\N}{\mathbb{N}}

\newcommand{\pt}{\partial_t}\newcommand{\pa}{\partial}
\newcommand{\les}{{\lesssim}}

\newcommand{\beeq}{\begin{equation}}\newcommand{\eneq}{\end{equation}}

\theoremstyle{plain}
\newtheorem{thm}{Theorem}[section]

\newtheorem{lem}[thm]{Lemma}

\theoremstyle{remark}
\newtheorem{rem}{Remark}

\theoremstyle{definition}

\newenvironment{prf}{\noindent {\bf Proof.} }{\endprf\par}
\def \endprf{\hfill  {\vrule height6pt width6pt depth0pt}\medskip}

\numberwithin{equation}{section}

\begin{document}

\title[Almost GWP for SLW with almost critical regularity]
{Almost global existence for some semilinear wave equations with
almost critical regularity}

\author{ Daoyuan Fang}
\address{Department of Mathematics, Zhejiang University, Hangzhou,
310027, China} \email{dyf@zju.edu.cn}
\thanks{The authors were supported in part by NSFC 10871175 and
10911120383.}

\author{Chengbo Wang}
\address{Department of Mathematics, Johns Hopkins University, Baltimore,
Maryland 21218} \email{wangcbo@jhu.edu}

\date{}
\dedicatory{} \commby{}

\begin{abstract}For any subcritical index of regularity $s>3/2$,
we prove the almost global well posedness for the 2-dimensional
semilinear wave equation with the cubic nonlinearity in the
derivatives, when the initial data are small in the Sobolev space
$H^s\times H^{s-1}$ with certain angular regularity. The lifespan is
known to be sharp in general. The main new ingredient in the proof
is an endpoint version of the generalized Strichartz estimates in
the space $L^2_t L_{|x|}^\infty L^2_\theta ([0,T]\times \R^2)$. In
the last section, we also consider the general semilinear wave
equations with the spatial dimension $n\ge 2$ and the order of
nonlinearity $p\ge 3$.
\end{abstract}

\maketitle


\section{Introduction}
The purpose of this paper is, for the 2-dimensional semilinear wave
equation with the cubic nonlinearity in the derivatives, to prove
the almost global well posedness with low regularity and sharp
lifespan. The main new ingredient in the proof will be an endpoint
version of the generalized Strichartz estimates in the space $L^2_t
L_{|x|}^\infty L^2_\theta ([0,T]\times \R^2)$. As complement, we
will also consider the general semilinear wave equations with the
spatial dimension $n\ge 2$ and the order of nonlinearity $p\ge 3$ in
the last section.

Let $\Box\equiv \pt^2-\Delta$, $\pa=(\pt,\pa_x)$ and $P_\al$ be
polynomials for $\al\in \N^3$, we consider the following Cauchy
problem \beeq\label{SLW} \Box u =\sum_{|\al|=3} P_\al(u) (\pa u)^\al
\eneq on $[0,T]\times\R^{2}$, together with the initial data at time
$t=0$ \beeq\label{data} u(0,x)=u_0(x), \ \pt u(0,x) = u_1(x)\ .
\eneq

In the case of classical $C_0^\infty$ initial data with size of
order $\ep$, the almost global existence with
\beeq\label{lifespan}T_\ep \ge \exp(c\epsilon^{-2})\eneq (for some
small constant $c>0$) can be proved by the standard energy methods,
see e.g. Sogge \cite{So08}. Moreover, the lifespan $T_\ep$ is also
sharp for the problem with nonlinearity $|\pt u|^3$ (Zhou
\cite{Zh01}).

Our object here is to prove the corresponding result with low
regularity. Note that the equation \eqref{SLW} with $P_\al$ be
constants $C_\al$ is invariant under the scaling transformation
$u(t,x)\rightarrow \la^{-\frac{1}{2}} u(\la t, \la x)$. This scaling
preserves the critical Sobolev space $\dot H^{s_c}$ with exponent
\beeq\label{sc}s_c = \frac{3}{2}\ ,\eneq which is then,
heuristically, a lower bound for the range of admissible $s$ such
that the problem \eqref{SLW}-\eqref{data} is well-posed in $C_t
H^s_x\cap C^1_t H^{s-1}_x$. (See e.g. Theorem 2 in \cite{FaWa08} for
the ill posed result with $s<s_c$ and nonlinearity $(\pt u)^3$.)

The local well posedness for the problem of this type with low
regularity has been extensively studied (see Ponce-Sideris
\cite{PoSi93}, Tataru \cite{Ta99} and the authors \cite{FaWa05}).
For this problem, besides scaling, there is one more mechanism due
to Lorentz invariance such that the problem is not well posed in
$C_t H^s_x\cap C^1_t H^{s-1}_x$ with $s=s_c+\ep$ for arbitrary small
$\ep\ll 1$ (see e.g. Lindblad \cite{Ld93}). Instead, the local well
posedness is true for $s>\frac{7}4$.

To state our main result, we need to introduce the Sobolev space
with angular regularity $b>0$, \beeq\label{ang-Sob}f\in
H^{s,b}_\theta \Leftrightarrow f\in H^s,\ \textrm{and
}(1-\pa_\theta^2)^{b/2} f\in H^s\eneq where the $(r,\theta)$ is the
polar coordinates. Now we are ready to state our main result.

\begin{thm}\label{thm-almost}
Let $n=2$, $s>s_c=3/2$ and $b>1/2$. Then there exist two small
positive constant $\ep_0$ and $c$, such that the problem
\eqref{SLW}-\eqref{data} admits an unique almost global solution
$(u,\pt u)\in C_{T_\ep}( H^{s,b}_\theta\times H^{s-1,b}_\theta)$
with $\pa_{t,x} u\in L^2_{T_\ep} L^\infty_{|x|} H^{b}_\theta$ on
$[0,T_\ep]\times \R^2$ with $T_\ep=\exp(c \ep^{-2})$, whenever
 $(u_0,u_1)\in
H^{s,b}_\theta\times H^{s-1,b}_\theta$ with norm bounded by $\ep\le
\ep_0$.
\end{thm}
\begin{rem}
  Here we note that, by adding some angular regularity, the Sobolev
  regularity required to ensure almost global existence
  is only $s>s_c$, which is $1/4$ less than the usual requirement of
  $s>7/4$.
\end{rem}

To prove Theorem~\ref{thm-almost} we shall require certain
Strichartz estimates, which involve the angular mixed-norm spaces
$$\|u\|_{L^q_t L^\infty_{|x|}L^2_\theta({\mathbb R}^2)}=
\left(\int_\R \textrm{esssup}_{\rho>0} \Bigl( \, \int_0^{2\pi}
|f(\rho(\cos \theta, \sin \theta))|^2 \, d\theta \, \Bigr)^{q/2}
dt\right)^{1/q}.$$

\begin{thm}\label{thm-Strichartz}  Let $P=\sqrt{-\Delta}$ in ${\mathbb
R}^2$. Then for any $\gamma>\frac12$, there exists a constant
$C_\gamma$ such that
\begin{equation}\label{eq-Stri-end}
\bigl\|\, e^{-itP}f \,  \bigr\|_{L^2_t L^\infty_{|x|}
L^2_\theta([0,T]\times {\mathbb R}^2)} \le C_\gamma
\left(\ln(2+T)\right)^\frac12\|f\|_{ H^{\gamma}({\mathbb R}^2)}.
\end{equation}
Moreover, if $2<q<\infty$, then
\begin{equation}\label{eq-Stri}
\bigl\|\, e^{-itP}f \,  \bigr\|_{L^q_t L^\infty_{|x|}
L^2_\theta({\mathbb R}\times {\mathbb R}^2)} \le C_{q}\|f\|_{\dot
H^{\gamma}({\mathbb R}^2)}, \quad \gamma= 1-1/q.
\end{equation}
\end{thm}

The estimate \eqref{eq-Stri-end} can be viewed as the endpoint
estimate of the estimates \eqref{eq-Stri}, which were proved
recently in Smith, Sogge and the second author \cite{SmSoWa09}. We
mention that the related estimates for \eqref{eq-Stri} where
$L^2_\theta$ is replaced by $L^r_\theta$ (with norms of different
regularity in the right) were proved by Sterbenz \cite{St05} for
$n\ge 4$ and the authors \cite{FaWa10} for the general case $n\ge
2$. In the radial case, the estimates \eqref{eq-Stri} and the higher
dimensional version were proved by the authors in \cite{FaWa06}
(with previous works in Sogge \cite{So08} for $n=3$, Sterbenz
\cite{St05} for $n\ge 3$).

\begin{rem} When $4<q<\infty$,
the Strichartz estimates \eqref{eq-Stri} is weaker than the standard
Strichartz estimates (see Theorem 3 of \cite{FaWa06})
\beeq\label{eq-Stri-cls} \bigl\|\, e^{-itP}f \, \bigr\|_{L^q_t
L^\infty_{x}({\mathbb R}\times {\mathbb R}^2)} \le C_{q}\|f\|_{\dot
H^{\gamma}({\mathbb R}^2)}, \quad \gamma= 1-1/q,\ 4<q<\infty. \eneq
By interpolating \eqref{eq-Stri} with \eqref{eq-Stri-cls}, we can
also improve $L^2_\theta$ to $L^p_\theta$ in \eqref{eq-Stri}.
\end{rem}

\begin{rem}\label{rem-3}
 Note that we have also the trivial energy estimate
\begin{equation}\label{energy}
\bigl\|\, e^{-itP}f\, \bigr\|_{L^\infty_t L^2_{|x|}L^2_\theta}\le
C\|f\|_{L^2},
\end{equation}
since $e^{-itP}$ is an unitary operator on $L^2$. By interpolation,
we can also get more general Strichartz type estimates involving
$L^q_t L^r_{|x|}L^2_\theta$ norm, where
$$\|u\|_{L^q_t L^r_{|x|}L^2_\theta(\R\times {\mathbb R}^2)}=
\left(\int_\R \left(\, \int_0^\infty \, \Bigl( \, \int_0^{2\pi}
|f(\rho(\cos \theta, \sin \theta))|^2 \, d\theta \, \Bigr)^{r/2}\,
\rho d\rho\, \right)^{q/r} d t\right)^{1/q}.$$
\end{rem}

This paper is organized as follows. In the next section, we give a
proof of Theorem~\ref{thm-Strichartz}, inspired by the arguments of
Smith, Sogge and the second author \cite{SmSoWa09} and Sterbenz
\cite{St05}. In Section 3, we shall prove the almost global
wellposedness with small data for the problem
\eqref{SLW}-\eqref{data}. To deal with the general nonconstant
functions $P_\al$, we will need to obtain improved bound for the
solution $u$, which is achieved in Lemma \ref{thm-bound-improved}.
In the final section, for the semilinear wave equations with general
dimension and general nonlinearity, we shall exploit further the
applications of the classical Strichartz estimates and their angular
improvement in Theorem~\ref{thm-Strichartz}, as an appendix to the
wellposed result for the 2-dimensional cubic semilinear wave
equation.

\section{Strichartz estimates}
In this section, we prove Theorem \ref{thm-Strichartz}, including
the critical $L^2_t L^\infty_{|x|} L^2_\theta$ Strichartz estimates
for the wave equation when $n=2$. We split the proof into three
steps. Although these steps are essentially the same except the last
step as in \cite{SmSoWa09}, we write out the complete proof for the
sake of completeness.

\subsection{Frequency Localization}
At first, we want to reduce the inequalities to the frequency
localized counterparts.

It is easy to see that the frequency localized estimates for Theorem
\ref{thm-Strichartz} are as follows
\begin{equation}\label{a}
\|e^{-itP}f\|_{L^q_tL^\infty_{|x|}L^2_\theta({\mathbb R}\times
{\mathbb R}^2)} \le C_q \|f\|_{L^2({\mathbb R}^2)}, \, \, \text{if }
\, q>2, \, \, \text{and } \, \Hat f(\xi)=0, \, |\xi|\notin [1/2,1]
\end{equation}
and \begin{equation}\label{a-end}
\|e^{-itP}f\|_{L^2_tL^\infty_{|x|}L^2_\theta([0,T]\times {\mathbb
R}^2)} \le C(\ln(2+T))^{1/2} \|f\|_{L^2({\mathbb R}^2)}, \, \,
\text{if } \, \Hat f(\xi)=0, \, |\xi|\notin [1/2,1]
\end{equation}
if $P=\sqrt{-\Delta}$.

By scaling and Littlewood-Paley theory, we see that \eqref{eq-Stri}
and \eqref{a} are equivalent. To deduce \eqref{eq-Stri-end} from
\eqref{a-end}, we will need to verify the following estimate for any
$\delta>0$ \beeq\label{eq-obs}\sum_{j\in \Z} 2^{j/2}
(1+2^j)^{-1/2-\delta} (\ln(2+2^j T))^{1/2}\le C_\delta
(\ln(2+T))^{1/2}.\eneq

In fact, if $T\ge e$, we deal with the following two different
cases.\\
 i) $2^j\ge 1$;
\begin{eqnarray*}
\sum_{j\ge 0} 2^{j/2}
(1+2^j)^{-1/2-\delta} (\ln(2+2^j T))^{1/2} &\le& \sum_{j\ge 0} 2^{-j \delta} (\ln(2+2^j T))^{1/2}\\
&\le& C
 \sum_{j\ge 0} 2^{-j\delta} (j\ln 2+\ln T)^{1/2}\\
 &\le& C \sum_{j\ge 0} 2^{-j\delta} (j\ln 2+1)^{1/2}(\ln T)^{1/2}\\
&\le& C_\delta (\ln T)^{1/2}.
\end{eqnarray*}
ii) $2^j\le 1$;
\begin{eqnarray*}
\sum_{j< 0} 2^{j/2}
(1+2^j)^{-1/2-\delta} (\ln(2+2^j T))^{1/2} &\le& \sum_{j< 0} 2^{j/2} (\ln(2+2^j T))^{1/2}\\
&\le&
 \sum_{j< 0} 2^{j/2} (\ln (2+ T))^{1/2}\\
 &\le& C (\ln (2+ T))^{1/2}.
\end{eqnarray*}

Else, if $T\le e$, we also deal with two different cases.\\
i) $2^j\ge T^{-1}$;
\begin{eqnarray*}
\sum_{2^j T\ge 1} 2^{j/2}
(1+2^j)^{-1/2-\delta} (\ln(2+2^j T))^{1/2} &\le& \sum_{2^j T\ge 1} 2^{-j \delta} (\ln(2+2^j T))^{1/2}\\
&\le& T^\delta\sum_{2^j T\ge 1} (2^j T)^{- \delta} (\ln(2+2^j T))^{1/2}\\
 &\le& C_\delta T^\delta\le \tilde{C}_\delta.
\end{eqnarray*}
ii) $1\le \la =2^j\le T^{-1}$;
$$\sum_{2^j T<1} 2^{j/2}
(1+2^j)^{-1/2-\delta} (\ln(2+2^j T))^{1/2} \le C \sum_{j} 2^{j/2}
(1+2^j)^{-1/2-\delta}\le C_\delta.$$

\subsection{Further reduction} Let us turn
to the proof of \eqref{a} and \eqref{a-end}.  Due to the support
assumptions for $\Hat f$ we have that
\begin{equation}\label{c}
\|f\|^2_{L^2({\mathbb R}^2)}\approx \int_0^\infty \int_0^{2\pi}|\Hat
f(\rho(\cos\omega,\sin\omega))|^2 d\omega d\rho.
\end{equation}
If we expand the angular part of $\Hat f$ using Fourier series we
find that if $\xi = \rho(\cos\omega,\sin\omega)$ then there are
Fourier coefficients $c_k(\rho)$ which vanish when $\rho\notin
[1/2,1]$ so that
$$\Hat f(\xi)=\sum_k c_k(\rho)e^{ik\omega},$$
and so, by \eqref{c} and Plancherel's theorem for $\mathbb{S}^1$ and
${\mathbb R}$ we have
\begin{equation}\label{d}
\|f\|^2_{L^2({\mathbb R}^2)} \approx \sum_k\int_{\mathbb
R}|c_k(\rho)|^2\, d\rho \approx \sum_k\int_{{\mathbb R}} |\Hat
c_k(s)|^2 \, ds,
\end{equation}
if $\Hat c_k$ denotes the one-dimensional Fourier transform of
$c_k(\rho)$. Recall that (see Stein and Weiss \cite{StWe71} p. 137)
\begin{equation}\label{e}
f(r(\cos\theta,\sin\theta))=\frac{1}{2\pi}\sum_k  \Bigl(\,
i^k\int_0^\infty J_k(r\rho)c_k(\rho)\rho \, d\rho\, \Bigr)
e^{ik\theta},
\end{equation}
if $J_k$ is the $k$-th Bessel function, i.e.,
\begin{equation}\label{f}
J_k(y)=\frac{(-i)^k}{2\pi}\int_0^{2\pi} e^{iy\cos\theta -ik\theta}\,
d\theta.
\end{equation}

Because of \eqref{e} and the support properties of the $c_k$, we
find that if we fix $\beta\in C^\infty_0({\mathbb R})$ satisfying
$\beta(\tau)=1$ for $1/2\le \tau\le 1$ but $\beta(\tau)=0$ if
$\tau\notin [1/4,2]$ then if we set
$\alpha=\rho\beta(\rho)\in{\mathcal S}({\mathbb R})$, we have
\begin{align*}
(e^{-itP}f)&(r(\cos\theta,\sin\theta))
\\
&= \frac{1}{2\pi}\sum_k \Bigl(\, i^k \int_0^\infty
J_k(r\rho)e^{-it\rho}c_k(\rho)\beta(\rho)\rho\, d\rho\,
\Bigr)e^{ik\theta}
\\
&=\frac1{(2\pi)^2}\sum_k\Bigl(\, i^k\int_0^\infty
\int_{-\infty}^\infty J_k(r\rho)e^{i\rho(s-t)} \Hat
c_k(s)\alpha(\rho)\, ds d\rho\, \Bigr)e^{ik\theta}
\\
&=\frac1{(2\pi)^3}\sum_k \Bigl(\, \int_0^\infty\int_{-\infty}^\infty
\int_0^{2\pi} e^{i\rho
r\cos\vartheta}e^{-ik\vartheta}e^{i\rho(s-t)}\Hat c_k(s)
\alpha(\rho)\, d\vartheta ds d\rho\, \Bigr)e^{ik\theta}
\\
&=\frac1{(2\pi)^3}\sum_k\Bigl(\,  \int_{-\infty}^\infty
\int_0^{2\pi} e^{-ik\vartheta}\Hat\alpha\bigl(\,
(t-s)-r\cos\vartheta \, \bigr)\, \Hat c_k(s) \, d\vartheta ds \,
\Bigr)e^{ik\theta}
\\
&=\frac1{(2\pi)^3}\sum_k\Bigl(\,  \int_{-\infty}^\infty \Hat c_k(s)
\psi_k(t-s,r) ds \, \Bigr)e^{ik\theta},
\end{align*}
where we set \beeq\label{eq-def}\psi_k(m,r)=\int_0^{2\pi}
e^{-ik\theta}\Hat\alpha\bigl(\, m-r\cos\theta \, \bigr)\,  \,
d\theta .\eneq

As a result, we have that for any $r\ge0$ \beeq\label{g}
\int_0^{2\pi}\Bigl| \, (e^{-itP}f)(r(\cos\theta,\sin\theta))\,
\Bigr|^2 \, d\theta =\frac1{(2\pi)^5} \sum_k\Bigl|\,
\int_{-\infty}^\infty \Hat c_k(s) \psi_k(t-s,r) \, ds\, \Bigr|^2.
\eneq

Now we claim that we have the estimate
\beeq\label{KeyStep}\|\psi_k(m,r)
\left<m\right>^{\frac12}\|_{L_m^2}\le C\, ,\eneq where
$\<m\>=\sqrt{1+m^2}$ and $C$ is independent of $k\in \Z$ and $r\ge
0$. If this is true, then
\begin{eqnarray*}
\|(e^{-itP}f)(r,\theta)\|_{L^2_\theta}&\le& C
\| \Hat c_k(s) \psi_k(t-s,r) \|_{l_k^2 L^1_s}\\
&\le& C \|\Hat c_k(s) \left<t-s\right>^{-1/2}\|_{l_k^2
L^2_s}\|\left<t-s\right>^{\frac12} \psi_k(t-s,r) \|_{L_s^2}\\
&\le& C \|\Hat c_k(s) \left<t-s\right>^{-1/2}\|_{l_k^2 L^2_s},
\end{eqnarray*}
and we can immediately get the required estimates \eqref{a} and
\eqref{a-end}, if we note that $$\left<t-s\right>^{-1/2}\in
L^q\textrm{ if }q>2\textrm{, and
}\|\left<t-s\right>^{-1/2}\|_{L^2_{t\in [0,T]}}\le C
(\ln(2+T))^{1/2}.$$

\subsection{The estimate for $\psi_k(m,r)$} Now we present the proof
of the key estimate \eqref{KeyStep} for $\psi_k(m,r)$, to conclude
the proof of the Strichartz estimate in Theorem
\ref{thm-Strichartz}.

We begin with the proof of the following pointwise estimates (which
is precisely Lemma 2.1 of \cite{SmSoWa09}).
\begin{lem}\label{mainest}  Let $\alpha\in {\mathcal S}({\mathbb R})$
and $N\in {\mathbb N}$ be fixed.
Then there is a uniform constant $C_N$, which is independent of
$m\in {\mathbb R}$ and $r\ge0$ so that the following inequalities
hold. First,
\begin{equation}\label{1}
\int_0^{2\pi} |\Hat\alpha(m-r\cos\theta)|\, d\theta \le C_N
\<m\>^{-N}, \quad \text{if } \, \, 0\le r\le 1, \, \, \text{or }\,
\, |m|\ge 2r.
\end{equation}
If $r>1$ and $|m|\le 2r$ then
\begin{equation}\label{2}
\int_0^{2\pi}|\Hat\alpha(m-r\cos\theta)|\, d\theta \le C\Bigl(\,
r^{-1}+r^{-1/2}\langle \, r-|m|\, \rangle^{-1/2}\, \Bigr).
\end{equation}
Consequently, for any $\delta>0$, we have the weaker estimate for
\eqref{KeyStep} \beeq\label{KeyStep-weaker} \|\psi_k(m,r)
\left<m\right>^{\frac12-\delta}\|_{L_m^2}\le C_\delta\, ,\eneq with
the constant $C_\delta$ independent of $r>0$.
\end{lem}

\begin{prf}
We first realize that \eqref{1} is trivial since $\Hat\alpha\in
{\mathcal S}$.  To prove \eqref{2}, it suffices to show that
\begin{equation}\label{4}
\int_0^{\pi/4}|\Hat\alpha(m-r\cos\theta)|\, d\theta
+\int_{\pi-\pi/4}^\pi |\Hat\alpha(m-r\cos\theta)|\, d\theta \le C
r^{-1/2}\langle \, r-|m|\, \rangle^{-1/2},
\end{equation}
and also
\begin{equation}\label{5}
\int_{\pi/4}^{\pi-\pi/4}|\Hat\alpha(m-r\cos\theta)| \, d\theta \le
Cr^{-1}.
\end{equation}

In order to prove \eqref{4}, it suffices to prove that the first
integral is controlled by the right side.  For if we apply this
estimate to the function $\Hat\alpha(-s)$, we then see that the
second integral satisfies the same bounds.  We can estimate the
first integral if we make the substitution $u=1-\cos\theta$, in
which case, we see that it equals
\begin{align*}
\int_0^{1-1/\sqrt2}|\Hat\alpha((m-r)+ru)|\, \frac{du}{\sqrt{2u-u^2}}
&\le \int_0^{1-1/\sqrt2}|\Hat\alpha((m-r)+ru)|\, \frac{du}{\sqrt u}
\\
&\le Cr^{-1/2}\int_0^\infty |\Hat\alpha((m-r)+u)|\, \frac{du}{\sqrt
u}
\\
&\le C'r^{-1/2}\, \langle \, r-m \, \rangle^{-1/2}
\\
&\le C'r^{-1/2}\, \langle \, r-|m| \, \rangle^{-1/2},
\end{align*}
as desired, which completes the proof of \eqref{4}.

To prove \eqref{5} we just make the change of variables
$u=r\cos\theta$ and note that $|du/d\theta| \approx r$ on the region
of integration, which leads to the inequality as $\Hat\alpha\in
{\mathcal S}$.

Finally, we check that inequalities \eqref{1} and \eqref{2} imply
\eqref{KeyStep-weaker}. If $r\le 1$, it is trivial. Else, if $r\ge
1$, we can prove \eqref{KeyStep-weaker} as follows
\begin{eqnarray*}
\|\psi_k(m,r) \left<m\right>^{\frac12-\delta}\|_{L_m^2}^2&\le&
 C+ C
\int_{|m|\le 2r} r^{-2} \left<m\right>^{1-2\delta} dm \\
&&+ C \int_{|m|\le 2r} r^{-1}  \left<m\right>^{1-2\delta}
\left<r-|m|\right>^{-1}  dm\\
&\le&
 C+ C
r^{-2} \left<r\right>^{2-2\delta}\\
&&+ C \int_{r/2\le m\le 2r} r^{-1} \left<m\right>^{1-2\delta}
\left<r-m\right>^{-1}  dm\\
&&+C \int_{0\le m\le r/2} r^{-1} \left<m\right>^{1-2\delta}
\left<r-m\right>^{-1}  dm\\
&\le& C+C \int_{r/2\le m\le 2r} r^{-2\delta} \left<r-m\right>^{-1} d
m\\&&+ C\int_{0\le m\le r/2} r^{-2} \left<m\right>^{1-2\delta} d m
\\&\le& C+C r^{-2\delta}\ln (2+r)\le C_\delta \ (\textrm{if } \delta>0).
\end{eqnarray*}

\end{prf}

Here, we remark that the reason we need to introduce a parameter
$\delta>0$ is due to the estimate \eqref{4} (the bound $r^{-1}$ will
be enough for us to get the estimate with $\delta=0$).

To prove the stronger estimate \eqref{KeyStep}, we need to consider
the effect of oscillated factor $e^{-ik\theta}$ in the definition of
$\psi_k(m,r)$, and the support property of the function $\alpha$.

To begin, we give some more reductions. At first, without loss of
generality, we can assume $m\ge 0$. In this case, we need only to
give the estimate for $\theta\in[0,\frac{3\pi}4]$ and
$\theta\in[\frac{3\pi}4,\pi]$. For the case $\theta\in
[\frac{3\pi}4,\pi]$, since $m-r\cos\theta\simeq m+r$ and $\Hat\al\in
\mathcal{S}$, the estimate is admissible for our purpose. So we need
only to give a refined estimate for the integral of the type
\begin{equation}\label{crucialterm}I_k(m,r)=\int_0^{3\pi/4} e^{-ik\theta}\Hat \alpha\bigl(\,
m-r\cos\theta \, \bigr)\,  \, d\theta \end{equation} when $m\le 2r$,
$r>1$. Moreover, we observe from \eqref{4} and \eqref{5} that
$$|\psi_k(m,r)|,\ |I_k(m,r)|\le C r^{-1},\textrm{ if }|m|\le r/2,$$ which are also
admissible estimates. This means that we need only to consider the
case $r/2<m<2r$ with $r>1$. Now we are ready to give the second
estimate about $\psi_k(m,r)$ (which is resemble to Proposition 4.1
of \cite{St05}).

\begin{lem}\label{mainest2}
Let $\al\in \mathcal{S}$ with support in $[1/4,2]$, $r>1$ and
$r/2<m<2r$. If $r<m+1$, then for any $N\ge 0$, we have
\begin{equation}\label{6}|I_k(m,r)|\le C_N r^{-1/2}
\left<r-m\right>^{-N},\end{equation} Else if $r>m+1$ and set
$d=\sqrt{r^2-m^2}$, we have \beeq\label{7}|I_k(m,r)|\le C r^{-1/2}
\left<r-m\right>^{-1/2}\left(\left<r-m\right>^{-1}+\min(k/d,
d/k)\right).\eneq Here, when $k=0$, the estimate is understood to be
$|I_0(m,r)|\le C r^{-1/2} \left<r-m\right>^{-3/2}$.
\end{lem}

Before giving the proof of Lemma~\ref{mainest2}, we give the proof
of \eqref{KeyStep}. By Lemma~\ref{mainest}, Lemma~\ref{mainest2} and
the discussion before Lemma~\ref{mainest2}, we know that
$$|\psi_k(m,r)|\le C\left\{\begin{array}{ll}
\left<m\right>^{-N}& |m|\ge 2r \textrm{ or } r\le 1\\
r^{-1}& |m|\le r/2 \textrm{ and } r> 1\\
\left<r+|m|\right>^{-N}+r^{-1/2}\left<|m|-r\right>^{-N} &
r<|m|+1, r/2\le|m|\le 2r \textrm{ and } r> 1\\
\left<r+|m|\right>^{-N}+r^{-1/2}\left<|m|-r\right>^{-3/2} &\\
\, \, \, \, +r^{-1/2}\left<|m|-r\right>^{-1/2}\min (k/d, d/k)& r\ge
|m|+1, r/2\le|m|\le 2r \textrm{ and } r> 1
\end{array}\right.
$$ Then a simple calculation will give us the key estimate \eqref{KeyStep}.

In fact, the case when $r\le 1$ is trivial. So we need only to
consider the case with $r>1$, in which case, we write the integral
into the sums as follows
\begin{eqnarray*}&& \int_\R |\psi_k(m,r)|^2
\left<m\right> d m\\
&=&\left(\int_{|m|\le r/2}+ \int_{|m|\ge 2r}+
\int_{\max(r-1,r/2)<|m|<2r}+  \int_{r/2<|m|<r-1}\right)
|\psi_k(m,r)|^2 \left<m\right> d m \\
&=&I+II+III+IV\ .\end{eqnarray*} The first two terms $I$ and $II$
can be estimated as before. For $III$,
 \begin{eqnarray*}
   III&\le & C+ \int_{\max(r-1,r/2)<|m|<2r}
   r^{-1}\left<|m|-r\right>^{-2N} \left<m\right> d m\\
   &\le & C + C \int_{\max(r-1,r/2)<|m|<2r}
   \left<|m|-r\right>^{-2N}  d m\le C.
 \end{eqnarray*}
Now we turn to the estimate for $IV$,
\begin{eqnarray*}
   IV&\le & C+ \int_{r/2<|m|<r-1}
   (r^{-1}\left<|m|-r\right>^{-3}+r^{-1}\left<|m|-r\right>^{-1}
   \min (k^2/d^2, d^2/k^2) )\left<m\right> d m\\
   &\le & C + C \int_{r/2<|m|<r-1}
   \left<|m|-r\right>^{-3}  d m \\
   &&+C \int_{r/2<|m|<r-1} \left<|m|-r\right>^{-1}\min (k^2/d^2, d^2/k^2) d m  \\
   &\le& C,
 \end{eqnarray*}
where in the last inequality, we used the fact that
$$\int_{r/2<|m|<r-1} \left<|m|-r\right>^{-1}\min (k^2/d^2, d^2/k^2) d
m\les 1.$$ In fact, if $k^2\le r$, then
\begin{eqnarray*}  && \int_{r/2<|m|<r-1}
\left<|m|-r\right>^{-1}\min (k^2/d^2, d^2/k^2) d m\\
& \les &\int_{r/2<|m|<r-1}
\left<|m|-r\right>^{-1}k^2 d^{-2} d m \\
&\le& C \int_{r/2<|m|<r-1} \left<|m|-r\right>^{-2}k^2 r^{-1}
 d m \\
&\le& C k^2/r \le C.
\end{eqnarray*}
Else, if $k^2>r$, we have
\begin{eqnarray*}  && \int_{r/2<|m|<r-1}
\left<|m|-r\right>^{-1}\min (k^2/d^2, d^2/k^2) d m\\
& \le &\int_{|m|<r-k^2/r}
\left<|m|-r\right>^{-1}k^2 d^{-2} d m \\
&&+\int_{\max(r/2, r-k^2/r)<|m|<r-1}
\left<|m|-r\right>^{-1}d^2 k^{-2} d m \\
&\le& C \int_{|m|<r-k^2/r} \left<|m|-r\right>^{-2}k^2 r^{-1}
 d m \\
&&+ C\int_{\max(r/2, r-k^2/r)<|m|<r-1}
r k^{-2} d m \\
&\le& C \left<k^2/r\right>^{-1}k^2/r + C rk^{-2} \min(k^2/r-1,
r/2-1)\le \tilde C.
\end{eqnarray*}
This proves our key estimate \eqref{KeyStep}.

Finally, we give the proof of Lemma~\ref{mainest2}, which will
conclude the proof of \eqref{KeyStep}.

\begin{prf}
If $r<m+1$, we have
$$m-r \cos\theta=r(1-\cos\theta)+m-r\ge m-r\ge -1.$$
Let $u=1-\cos\theta$, so we get
$$\left<m-r \cos\theta\right>\simeq 1+r(1-\cos\theta)+(2+m-r)\simeq
\left<m-r\right>+\left<ru\right>.$$
Since $\alpha\in {\mathcal S}$, \begin{eqnarray*}   |I_k(m,r)|&\le&
\int_0^{\frac{3\pi}4}|\Hat \alpha\bigl(\, m-r\cos\theta \, \bigr)|\,
\, d\theta \\
&\le& C \int_0^{\frac{3\pi}4}
\left<m-r\cos\theta\right>^{-2N}\, d\theta \\
&\le& C \int_0^{\frac{3\pi}4}
\left<m-r\right>^{-N}\left<ru\right>^{-N}\, d\theta \\
&\le& C
\int_0^{1+1/\sqrt{2}}\left<m-r\right>^{-N}\left<ru\right>^{-N}
\, \frac{du}{\sqrt{2u-u^2}} \\
&\le& C
\int_0^{1+1/\sqrt{2}}\left<m-r\right>^{-N}\left<ru\right>^{-N}
\, \frac{du}{\sqrt{u}} \\
&\le& C r^{-1/2}\left<m-r\right>^{-N},
\end{eqnarray*} which gives us \eqref{6}.

Now we turn to the proof for the case $r\ge m+1$. We can imagine
that the behavior is worst in the region that $m-r\cos\theta\sim 0$.
To illustrate this, we introduce $\theta_0\in (0,\frac\pi2]$ such
that \beeq\label{eq-theta0-def}r\cos\theta_0=m,\
\sin\theta_0=\frac{\sqrt{r^2-m^2}}r\equiv \frac dr.\eneq Then the
local behavior of the function $r\cos\theta-m$ near
$\theta=\theta_0$ looks like
$$r\cos\theta-m\sim -d (\theta-\theta_0)+\mathcal{O}((\theta-\theta_0)^2),$$
since $r\cos\theta_0-m=0$ and
$\frac{d}{d\theta}(r\cos\theta-m)|_{\theta=\theta_0}=-r\sin\theta_0=-d$.
Based on this information, we make the change of variable
\beeq\label{eq-phi-defi}\be=d(\theta-\theta_0),\
\phi(\be)=m-r\cos(\theta_0+\be/d)\ .\eneq For the function $\phi$,
we can find that $\phi(0)=0$, $\phi'(\beta)=\frac rd \sin\theta$.
Moreover, we have the following
\begin{lem}\label{lem-phi-pro}
Let $\phi(\be)$ be the function defined by \eqref{eq-phi-defi}, and
$\theta\in [\theta_1, \frac{3\pi}4]$ with $\theta_1\in
(0,\frac{\pi}4)$ such that $\frac rd \sin\theta_1=\frac12$, then we
have \beeq\label{eq-phi}\frac12\le \phi'(\be)\le 1+|\phi(\be)|\simeq
\left<\phi(\be)\right>.\eneq In addition,
$$|\phi(\beta)|=|\phi(\be)-\phi(0)|\ge \frac12 |\be|.$$
\end{lem}
\begin{prf}
We need only to give the proof of the inequality
$$\phi'(\be)\le 1+|\phi(\be)|.$$
 In fact, if
$\be=0$ (i.e., $\theta=\theta_0$), we know the inequality is true
with identity (see \eqref{eq-theta0-def}). For $\be\le 0$, we have
$\phi(\be)\le 0$ and the inequality amounts to $\phi'(\be)\le
1-\phi(\be)$, which is equivalent to
$$\frac rd \sin\theta\le 1+r\cos\theta-m,\ \theta\in [\theta_1,\theta_0].$$
Now we can see that this inequality is trivial by the monotonicity
of the trigonometric functions
$$\frac rd \sin\theta\le \frac rd \sin\theta_0=1=1+r\cos\theta_0-m \le 1+r\cos\theta-m.$$
 If we consider instead the case $\be\ge 0$, we know that it is
 equivalent to \beeq\label{eq-phi-p}\frac rd \sin\theta\le 1+m-r\cos\theta,\ \theta\in
 [\theta_0,\frac{3\pi}4].\eneq
 Once again, by the monotonicity of the trigonometric functions, we
need only to prove the inequality for $\theta\in [\theta_0,
\frac{\pi}2]$. In the latter case, consider $F(\theta)=
1+m-r\cos\theta-\frac rd \sin\theta$, we observe that (recall $r\ge
m+1$)
$$F'(\theta)=r\sin\theta-\frac rd \cos\theta\ge r\sin\theta_0-\frac rd \cos\theta_0=d-\frac md
=\frac{d^2-m}d\ge \frac{(m+1)^2-m^2-m}d\ge 0\ .$$ Recall
$F(\theta_0)=0$, we know that $F(\theta)\ge 0$ for $\theta\in
[\theta_0, \frac{\pi}2]$ and so is \eqref{eq-phi-p}. This completes
the proof of the inequality \eqref{eq-phi}.
\end{prf}

Now let us continue the proof of the estimate for $I_k$. We write
\beeq\label{I=J+K} I_k(m,r)=J_k(m,r)+K_k(m,r)\eneq with
\begin{eqnarray}
J_k(m,r)&=&\int_{\theta_1}^{3\pi/4} e^{-ik\theta}\Hat \alpha\bigl(\,
m-r\cos\theta \, \bigr)\,  \, d\theta\nonumber\\
&=&\frac{e^{-ik\theta_0}}{d}\int_{d(\theta_1-\theta_0)}^{d(3\pi/4-\theta_0)}
e^{-i\frac kd\be}\Hat \alpha\bigl(\phi(\be)\bigr)\, \, d\be
\equiv\frac{e^{-ik\theta_0}}{d} L_k.\label{J-L}
\end{eqnarray}

We first give the easier estimate for $K_k$. In fact, if $\theta\in
[0,\theta_1]$,
$$r\cos\theta-m\ge r\cos\theta_1-m=\frac{\sqrt{3r^2+m^2}}2-m\ge
\frac{3(r^2-m^2)}{4\sqrt{3r^2+m^2}}\simeq r-m .$$ Note that
$\theta_1\sim \sin\theta_1=\frac{d}{2r}$, this means that
\beeq\label{K-k}|K_k(m,r)|\le C \int_0^{\theta_1}
(r-m)^{-N-1}d\theta\le C \frac dr (r-m)^{-N-1}\le C
r^{-1/2}(r-m)^{-N}.\eneq

Now we turn to the estimate for $J_k$ in terms of $L_k$. We want to
exploit the effect of oscillated factor $e^{-i\frac kd\be}$,
together with the support property of the function $\alpha$. Recall
that $i\frac dk\partial_\be e^{-i\frac kd\be}= e^{-i\frac kd\be}$,
we use integration by parts in $\be$ to get
\begin{eqnarray*}
  |L_k(m,r)|&=&\left|\frac dk\int_{d(\theta_1-\theta_0)}^{d(3\pi/4-\theta_0)}
\partial_\be (e^{-i\frac kd\be})\Hat \alpha\bigl(\phi(\be)\bigr)\, \, d\be\right|\\
&=&\frac dk \left|e^{-i\frac kd\be}\Hat
\alpha\bigl(\phi(\be)\bigr)|_{d(\theta_1-\theta_0)}^{d(3\pi/4-\theta_0)}-
\int_{d(\theta_1-\theta_0)}^{d(3\pi/4-\theta_0)} e^{-i\frac kd\be}
\phi'(\be) (\Hat \alpha)'\bigl(\phi(\be)\bigr)\,
\, d\be\right|\\
&\le &C \frac dk \left(|\Hat \alpha\bigl(m-r\cos \frac {3\pi}
4\bigr)|+|\Hat \alpha\bigl(m-r\cos \theta_1\bigr)|+\int_{\R}
\left<\phi(\be)\right>^{1-N}\,
\, d\be\right)\\
&\le &C \frac dk \left(1+\int_{\R} \left<\be\right>^{1-N}\,
\, d\be\right)\\
&\le &C \frac dk \ \textrm{, if }k\neq 0,
\end{eqnarray*} where we have used Lemma \ref{lem-phi-pro} in the
first and second inequality.

To prove another inequality for $|L_k|$, we need only to exploit the
support property of $\al$. Since $\textrm{supp} \al \subset [\frac
14, 2]$, we can introduce $\tilde \al (\rho)= i \al(\rho)/\rho\in
\mathcal{S}$ so that $\Hat \al=(\Hat{\tilde \al})'$ and
$$\Hat \al (\phi(\be))=(\Hat{\tilde \al})'(\phi(\be))=
\frac{1}{\phi'(\be)}\pa_\be(\Hat{\tilde \al}(\phi(\be))).
$$ Thus
we have
\begin{eqnarray*}
  |L_k(m,r)|&=&\left|\int_{d(\theta_1-\theta_0)}^{d(3\pi/4-\theta_0)}
e^{-i\frac kd\be} \frac{1}{\phi'(\be)}\pa_\be(\Hat{\tilde \al}(\phi(\be)))\, \, d\be\right|\\
&\le &\left|\left.\left(e^{-i\frac kd\be} \frac{1}{\phi'(\be)}
\Hat{\tilde
\al}(\phi(\be))\right)\right|_{d(\theta_1-\theta_0)}^{d(3\pi/4-\theta_0)}\right.\\
&&\left.- \int_{d(\theta_1-\theta_0)}^{d(3\pi/4-\theta_0)}
\pa_\be\left(e^{-i\frac kd\be}
\frac{1}{\phi'(\be)}\right)\Hat{\tilde \al}(\phi(\be))\,
\, d\be\right|\\
&\le &C \left(|\Hat{\tilde \al}\bigl(m-r\cos \frac{3\pi}
4\bigr)|+|\Hat{\tilde \al}\bigl(m-r\cos
\theta_1\bigr)|\right)\\
&&+\int_{d(\theta_1-\theta_0)}^{d(3\pi/4-\theta_0)}
\left|\pa_\be\left(e^{-i\frac kd\be}
\frac{1}{\phi'(\be)}\right)\Hat{\tilde \al}(\phi(\be))\,
\, \right|d\be\\
&\le &C \left<r-m\right>^{-N}+C\frac kd
\int_{d(\theta_1-\theta_0)}^{d(3\pi/4-\theta_0)} \left|
\frac{1}{\phi'(\be)}\Hat{\tilde \al}(\phi(\be))\,
\right|d\be\\
&&+C\int_{d(\theta_1-\theta_0)}^{d(3\pi/4-\theta_0)} \left|
\frac{\phi''(\be)}{(\phi'(\be))^2}\Hat{\tilde \al}(\phi(\be)) \,
\right|d\be\\
&\le &C\left<r-m\right>^{-1}+C\frac kd \ ,
\end{eqnarray*}where we have used the fact that
$r\cos\theta_1-m\gtrsim r-m$, $\phi'(\be)\ge \frac 12$ (for
$\theta\in [\theta_1, \frac{3\pi}4]$) and
$\phi''(\be)=\frac{r}{d^2}\cos\theta=\mathcal{O}((r-m)^{-1})$.
Combining with the previous inequality, we have proved
$$|L_k(m,r)|\le C\left<r-m\right>^{-1}+C\min\left(\frac kd,\frac dk\right)$$
and so is the inequality \eqref{7} by \eqref{I=J+K}, \eqref{J-L} and
\eqref{K-k}. This completes the proof.\end{prf}

\section{Almost global existence for cubic SLW} In this section, we
prove Theorem \ref{thm-almost}, as an application of the endpoint
estimate \eqref{eq-Stri-end}.

To begin, let us prove the fractional Leibniz rule in the Sobolev
space with angular regularity.
\begin{lem}\label{thm-Leibniz}
Let $n=2$, $s\in (0,1)$, $b>\frac{1}2$ and $\psi\in
\mathcal{S}(\R^2)$ be a radial function, then we have
\beeq\label{eq-infty} \|\psi*f\|_{L^\infty_r L^2_\theta}\les
\|\psi\|_{L^1}\|f\|_{L^\infty_r L^2_\theta}\ ,\eneq and the
fractional Leibniz rule
\beeq\label{eq-Leibniz}\|fg\|_{H^{s,b}_\theta}\les
\|f\|_{L^\infty_{|x|}H^b_\theta}\|g\|_{H^{s,b}_\theta}+
\|g\|_{L^\infty_{|x|}H^b_\theta}\|f\|_{H^{s,b}_\theta}\ .\eneq
Moreover, we have \beeq\label{eq-Leibniz-Riesz} \|fg\|_{
H^{s,b}_{\theta}}\les \|f\|_{ L^\infty_{|x|}H^b_\theta\cap \dot
H^{1,b}_{\theta}}\|g \|_{ H^{s,b}_\theta}\ ,\eneq
\beeq\label{eq-Leibniz-Riesz-1} \|fg\|_{\dot H^{1,b}_{\theta}\cap
L^\infty_{|x|}H^b_\theta}\les \|f\|_{\dot H^{1,b}_{\theta} \cap
L^\infty_{|x|}H^b_\theta}\|g \|_{\dot H^{1,b}_{\theta} \cap
L^\infty_{|x|}H^b_\theta}\ .\eneq
\end{lem}
\begin{prf}
At first, we give the proof for \eqref{eq-infty}. Recall $$(\psi*
f)(x)=\int \psi(y) f(x-y)dy,$$ we set $x=(r\cos
\omega,r\sin\omega)$, $y=(\la \cos \theta,\la \sin\theta)$, then
$x-y=(\rho \cos \al,\rho \sin\al)$, with
$$\rho=\sqrt{r^2+\la^2-2r\la \cos (\omega-\theta)},\ \al=\omega+\arcsin
\left(\frac{\la}{\rho}\sin (\omega-\theta)\right).$$
Introducing a new variable $a=\omega-\theta\in [0,2\pi]$, then
$\rho=\rho(\la, r, a)$ and $\al=\al(\la,
r,\omega,a)=\omega+h(\la,r,a)$ for some function $h$. Now, for fixed
$r$,
\begin{eqnarray*} \|\psi
*f\|_{L^2_\omega}&=&\|\int_{0}^\infty \int_0^{2\pi} \psi(\la) f(x-y)\la d\la
d\theta\|_{L^2_\omega}\\
& \le&  \| f(\rho\cos\al,\rho\sin\al) \|_{L^\infty_\la L^2_\omega
L^1_\theta} \int_{0}^\infty |\psi(\la)| \la d\la
\\
&\simeq&\|\psi\|_{L^1} \|f(\rho\cos\al,\rho\sin\al) \|_{L^\infty_\la
L^2_\omega
L^1_a}\\
&\les& \|\psi\|_{L^1} \|f(\rho(\la, r, a)\cos\al(\la, r,
\omega,a),\rho(\la, r, a)\sin\al(\la, r, \omega,a))
\|_{L^\infty_\la L^2_\omega L^2_a}\\
&\les& \|\psi\|_{L^1}\|f(\rho(\la, r, a)\cos\omega,\rho(\la, r, a)\sin\omega)
\|_{L^\infty_\la L^2_a L^2_\omega}\\
&\les& \|\psi\|_{L^1}\|f(\rho(\la, r, a)\cos\omega,\rho(\la, r, a)\sin\omega)
\|_{L^\infty_\la L^\infty_a L^2_\omega}\\
&\le&\|\psi\|_{L^1}\|f(\rho\cos\omega,\rho\sin\omega)
\|_{L^\infty_\rho L^2_\omega}\ ,\end{eqnarray*} which proves
\eqref{eq-infty}. The estimate \eqref{eq-infty} tell us that the
space $L^\infty_r L^2_\theta$ is stable under the frequency
localization.

Based on \eqref{eq-infty}, and the fact that
 $H^{b}_\theta$ is an algebra under multiplication when $b>1/2$, we
 can easily apply Littlewood-Paley decomposition to prove the
fraction Leibniz rule \eqref{eq-Leibniz}, \eqref{eq-Leibniz-Riesz},
and \eqref{eq-Leibniz-Riesz-1}.
\end{prf}

Now we are ready to give the proof of Theorem \ref{thm-almost},
based on the endpoint estimate \eqref{eq-Stri-end} and the fraction
Leibniz rule \eqref{eq-Leibniz}.

At first, we prove the easier case when $P_\al(u)$ do not depend on
$u$, for which the idea of the proof will be clear. After that, we
will modify the proof to show that we can still handle the general
case.

\subsection{The case with $P_\al(u)=C_\al$} By \eqref{eq-Stri-end} and energy estimate, for fixed
$s>\frac32$ and $b>\frac12$, we have \beeq\label{eq-Stri-u}\|e^{-i t
P} f\|_{L^2_T L^\infty_{|x|} H^{b}_\theta} \le C_0 (\ln(2+T))^{1/2}
\|f\|_{H^{s-1,b}_\theta}\eneq and \beeq\label{eq-energy} \|e^{-i t
P} f\|_{L^\infty_T H^{s-1, b}_\theta}\le C_0
\|f\|_{H^{s-1,b}_\theta}\eneq
 with some constant $C_0>1$.
Recall that we have the initial data $(u_0,u_1)\in
H^{s,b}_\theta\cap H^{s-1,b}_\theta$ with
\beeq\label{IData}\|u_0\|_{ H^{s,b}_\theta}+\|u_1\|_{
H^{s-1,b}_\theta}=\ep\le \ep_0,\eneq where $\ep_0$ will be fixed
later (see \eqref{cons-def}).

Given the metric \beeq\label{metric}d(u,v)=(\ln(2+T))^{-1/2}
\|\pa_{t,x}(u-v)\|_{L^2_T L^\infty_{|x|}
H^{b}_\theta}+\|\pa_{t,x}(u-v)\|_{ L^\infty_T H^{s-1,b}_\theta
},\eneq we define the complete domain with $T=\exp(c \ep^{-2})$ and
$c\ll 1$ to be chosen later (see \eqref{cons-def}), $$X=\{u\in C_T
H^{s,b}_\theta\cap C_T^1 H^{s-1,b}_\theta: u(0,x)=u_0(x),\ \pt
u(0,x)=u_1(x),\ d(u,0)<\infty\}. $$

Then for any $u\in X$, we denote $\Pi u$ to be the solution to the
linear wave equation
$$\Box \Pi u=(\pa u)^\al$$ with initial data $(u_0, u_1)$.
Note that for $u\in X$, we have $(\pa u)^\al\in L^1_T
H^{s-1,b}_\theta$ by using the fraction Leibniz rule
\eqref{eq-Leibniz}, and so $\Pi u\in C_T H^{s,b}_\theta\cap C_T^1
H^{s-1,b}_\theta$ is well defined. Thus, by energy estimates
\eqref{eq-energy} and Strichartz estimates \eqref{eq-Stri-u}, we
have
$$
  d(\Pi u, 0) \le C_1 (\|u_0\|_{H^{s,b}_\theta}+\|u_1\|_{H^{s-1,b}_\theta})+ C_1 \|(\pa u)^\al\|_{L^1_T
  H^{s-1,b}_\theta}$$ for some $C_1\ge C_0$.
Based on this estimate, we define a complete domain $D_\ep \subset
X$ so that the map $\Pi$ will be a contraction map in $D_\ep$ (for
$\ep_0$ and $c$ small enough),
$$D_\ep =\{u\in X; \ d(u,0)\le 2 C_1 \ep\}\ .$$

By using the fraction Leibniz rule \eqref{eq-Leibniz} and noting
that $T=\exp(c\ep^{-2})>e>2$, we have for some $C_2\ge C_1$,
\begin{eqnarray*}
  d(\Pi u, 0) &\le& C_1 (\|u_0\|_{H^{s,b}_\theta}+\|u_1\|_{H^{s-1,b}_\theta})+ C_1 \|(\pa u)^\al\|_{L^1_T H^{s-1,b}_\theta}\\
&\le & C_1 \ep+ C_2 \|\pa u\|_{L^2_T L^\infty_{|x|} H^{b}_\theta}^2  \|\pa u\|_{L^\infty_T H^{s-1,b}_\theta}\\
&\le & C_1 \ep+ C_2 \ln(2+T) d(u,0)^3\\
&\le & C_1 \ep+ C_2 d(u,0)^3 +C_2 c \ep^{-2} d(u,0)^3 \ .
\end{eqnarray*}
Moreover, for any $u,v\in X$, we have for some $C_3\ge C_2$,
\begin{eqnarray*}
  d(\Pi u, \Pi v) &\le& C_1 \|(\pa u)^\al-(\pa v)^\al\|_{L^1_T H^{s-1,b}_\theta}\\
&\le & C_3 \ln(2+T) (d( u,0)^2+ d( v,0)^2)  d(u,v)\\
&\le & (C_3+C_3  c \ep^{-2}) (d( u,0)^2+ d( v,0)^2) d(u,v).
\end{eqnarray*}

Now we fix the constant $\ep_0$ and $c$ such that we have
$$c\ep_0^{-2}\ge 1,\ 2 C_2 (2 C_1 \ep_0)^3 c \ep_0^{-2}\le C_1 \ep_0,\ 2 C_3 c \ep_0^{-2} (2\times 2 C_1 \ep_0)^2\le \frac 12,$$
which can be satisfied if we set \beeq\label{cons-def}c=\frac 1{2^6
C_1^2 C_3},\ \ep_0=\sqrt{c}\ .\eneq Then we know that if $\ep\le
\ep_0$,  the map $\Pi$ is a contraction map on the complete set
$D_\ep$, and
 the fixed point $u\in D_\ep$ of the map $\Pi$ gives the required
unique almost global solution.

\subsection{General case involving the unknown function}
In the general case, we need to give the estimate for $u$ in
$L^\infty_{t,x}$.

To begin, recall $u(t)=u(0)+\int_0^t \pt u (s)ds$, then we have the
trivial bound \beeq\label{eq-est-u-increase}\|u(t)\|_{L^2}\le
\|u(0)\|_{L^2}+t\|\pt u(s)\|_{L^\infty_{s\in [0,t]} L^2}\ .\eneq
This will give us the uniform bound in time of order $1$, if we
combine the bound of $\pa_{t,x} u$ in $L^\infty_t H^{s-1}$ ($s>1$).
Note that we are in the situation that for $t\le \exp(c \ep^{-2})$,
such a bound for $u$ will not be admissible for our purpose.

Instead, we need to prove an improved estimates on the $L^\infty$
norm on $u\in X$. We claim that we can in fact prove the uniform
bound for $u$, \beeq\label{eq-est-u}\|u\|_{L^\infty_{t\in
[0,T_0]}L^\infty_{|x|}H^b_\theta}\le C\eneq with $C$ independent of
$\ep$, if $T_0=\exp(\ep^{-2})$, $u\in D_\ep$ and
$\|u(0,\cdot)\|_{H^{s,b}_\theta}\le \ep$. In fact, this estimate is
a special case of the following lemma.
\begin{lem}\label{thm-bound-improved}
Let $n\ge 1$, $s>\frac n2$ and $s\ge 1$. Assume we have $u\in C
H^{s}\cap C^1 H^{s-1} ([0,T_0]\times \R^n)$,
 then there is an universal constant $C$ such
 that
\beeq\label{eq-est-infty-final}\|u(t)\|_{L^\infty}\le C \left(\frac
2n\ \frac{1-\delta}{\delta}\right)^{\frac{1-\delta}{2}}
(\|u(t)\|_{\dot H^{\frac{n}{2(1-\delta)}}}^{1-\delta}
\|u(0)\|_{L^2}^{\delta} + t^\delta \|\pa_{t,x} u\|_{L^\infty_t
H^{s-1}} ) \eneq for any $0<\delta\le 1-\frac n{2s}$ and $0\le t\le
T_0$. In particular, if $n=2$, $T_0=\exp(\ep^{-2})$ and
$\|u(0)\|_{L^2}+\|\pa_{t,x} u\|_{L^\infty_t H^{s-1}}\le C \ep$ with
$\ep\ll 1$, then by choosing $\delta=\ep^2$, we have a uniform bound
in $\ep$
$$\|u\|_{L^\infty_{t,x}}\les 1\ .$$ Moreover, if $n=2$, $T_0=\exp(\ep^{-2})$,
$u(0)=0$ and $\|\pa_{t,x} u\|_{L^\infty_t H^{s-1}}<\infty$, then by
choosing $\delta=\ep^2$, we have
\beeq\label{eq-est-u-conv}\|u\|_{L^\infty_{t,x}}\le C \ep^{-1}
\|\pa_{t,x} u\|_{L^\infty_t H^{s-1}}.\eneq
\end{lem}
\begin{prf}
Recall that we have ($0<\delta <1$ 
)
$$\|u\|_{L^\infty}\le C_\delta \|u\|_{\dot H^{\frac n{ 2(1-\delta)}}}^{1-\delta}\|u\|_{L^{2}}^{\delta}.$$
To prove the improved bound for $u$, we need to know the explicit
dependence of $C_\delta$ with respect to $\delta\in (0,1)$,
\begin{eqnarray*} \|f\|_{L^\infty}&\le&
C\|\hat f\|_{L^1}\\
&\le& C \||\xi|^{-\frac n2-a} |\xi|^{\frac
n2+a}\hat{f}(\xi)\chi_{|\xi|\ge
\la}\|_{L^1}+C \|\hat{f}(\xi)\chi_{|\xi|\le \la}\|_{L^1}\\
&\le& C \||\xi|^{-\frac n2-a}\chi_{|\xi|\ge \la}\|_{L^2}
\||\xi|^{\frac n2+a}\hat{f}(\xi)\|_{L^2}+C\la^{\frac n2} \|\hat{f}(\xi)\|_{L^2}\\
&\le& \frac{C}{\sqrt{a}}\la^{-a} \|f\|_{\dot H^{\frac
n2+a}}+C\la^{\frac n2} \|f\|_{L^2}\ .
\end{eqnarray*}
By choosing $$\la=\left(\frac{\|f\|_{\dot H^{\frac
n2+a}}}{\sqrt{a}\|f\|_{L^2}}\right)^{\frac1{\frac n2+a}},$$ we
conclude that
$$\|f\|_{L^\infty}\le C a^{-\frac{n}{2n+4a}}\|f\|_{\dot
H^{\frac n2+a}}^{\frac{\frac n2}{\frac n2+a}} \|f\|_{L^2}^{\frac
a{\frac n2+a}}$$ i.e. \beeq\label{eq-est-infty}\|f\|_{L^\infty}\le C
\left(\frac 2n\
\frac{1-\delta}{\delta}\right)^{\frac{1-\delta}{2}}\|f\|_{\dot
H^{\frac{n}{2(1-\delta)}}}^{1-\delta} \|f\|_{L^2}^{\delta}\eneq From
which we see that we can have the upper bound for $C_\delta$
$$C_\delta\le C \left(\frac 2n\ \frac{1-\delta}{\delta}\right)^{\frac{1-\delta}{2}}.$$

Now we can give the improved estimates for $u$. Observe that if we
choose $\delta>0$ small enough ($\delta\le 1-\frac n{2s}$ will be
enough), we will have $\frac n{2(1-\delta)}\le s$. Combining
\eqref{eq-est-infty} with \eqref{eq-est-u-increase}, we know that
$$\|u(t)\|_{L^\infty}\le C C_\delta (\|u(t)\|_{\dot
H^{\frac{n}{2(1-\delta)}}}^{1-\delta} \|u(0)\|_{L^2}^{\delta} +
t^\delta \|\pa_{t,x} u\|_{L^\infty_t H^{s-1}} ),
$$ which gives us \eqref{eq-est-infty-final}.
\end{prf}

Now we are ready to prove the general case of Theorem
\ref{thm-almost}. For any $u\in D_\ep$, we denote $\Pi u$ to be the
solution to the linear equation
$$\Box \Pi u=\sum_{|\al|=3}P_\al(u) (\pa u)^\al$$ with data $(u_0, u_1)$.
Recall that in the definition of the space $X$, we have assumed $T=
T_\ep \equiv \exp(c \ep^{-2}) \le \exp(\ep^{-2})$ since $c\le 1$,
and so that we have the bound \eqref{eq-est-u}
$$\|u\|_{L^\infty_{T}L^\infty_{|x|}H^b_\theta}\les 1 \ .$$ By
\eqref{eq-Leibniz-Riesz-1} and $\|u\|_{L^\infty_T \dot
H^{1,b}_{\theta}}\le d(u,0)\le 2 C_1 \ep \les 1$, we know that
$$\|P_\al(u)-P_\al(0)\|_{L^\infty_T \dot H^{1,b}_{\theta}\cap
L^\infty_{T} L^\infty_{|x|}H^b_\theta}\les 1.$$

Thus, by \eqref{eq-Leibniz} and \eqref{eq-Leibniz-Riesz}, together
with \eqref{eq-Stri-u} and \eqref{eq-energy}, we have for some
$\tilde C_2\ge C_2$
\begin{eqnarray*}
  d(\Pi u, 0) &\le& C_1 (\|u_0\|_{H^{s,b}_\theta}+\|u_1\|_{H^{s-1,b}_\theta})+
  C_1 \|P_\al(u) (\pa u)^\al\|_{L^1_T H^{s-1,b}_\theta}\\
  &\le& C_1 \ep+C_1 \| P_\al(0)(\pa u)^\al\|_{L^1_T H^{s-1,b}_\theta}+
  C_1 \|(P_\al(u)-P_\al(0)) (\pa u)^\al\|_{L^1_T H^{s-1,b}_\theta}\\
&\le & C_1 \ep+ 2 C_2 c \ep^{-2} d(u,0)^3 \max_\al |P_\al(0)|+
\tilde C_2 \|\pa u\|_{L^2_T L^\infty_{|x|} H^{b}_\theta}^2
\|\pa u\|_{L^\infty_T H^{s-1,b}_\theta}\\
&\le & C_1 \ep+ 4 \tilde C_2 c \ep^{-2} d(u,0)^3\ .
\end{eqnarray*}

Similarly, for any $u,v\in D_\ep$, we have $u(0)-v(0)=u_0-u_0=0$ and
so we can apply \eqref{eq-est-u-conv} for $u-v$ as follows
\beeq\label{eq-est-u-conv2}\|u-v\|_{L^\infty_{t,|x|}H^b_\theta}\le
\tilde C_4 \ep^{-1} \|\pa_{t,x}(u-v)\|_{L^\infty_t
H^{s-1,b}_\theta}.\eneq Recall again the fractional Leibniz rule
\eqref{eq-Leibniz}-\eqref{eq-Leibniz-Riesz-1}, we have for some
constant $\tilde C_3\ge C_3$,
\begin{eqnarray*}
  d(\Pi u, \Pi v) &\le& C_1 \|P_\al(u)(\pa u)^\al-P_\al(v)(\pa v)^\al\|_{L^1_T H^{s-1,b}_\theta}\\
&\le& C_1 \|P_\al(u)((\pa u)^\al-(\pa v)^\al)\|_{L^1_T
H^{s-1,b}_\theta}+C_1 \|(P_\al(u)-P_\al(v))(\pa v)^\al\|_{L^1_T
H^{s-1,b}_\theta}\\
&\le & \tilde C_3 \ln(2+T) (d( u,0)^2+ d( v,0)^2)  d(u,v)+ \tilde C_3 \ln(2+T) d( v,0)^3 \ep^{-1} d(u,v)\\
&\le & 4 \tilde C_3  c\ep^{-2} \left(d( u,0)^2+ d( v,0)^2+ \ep^{-1}
d( v,0)^3 \right) d(u,v).
\end{eqnarray*}

Now as in the case $P_\al=C_\al$, we can find some small positive
constants $c$ and $\ep_0$, such that the previous estimates amount
to
$$d(\Pi u,0)\le 2 C_1 \ep \textrm{ and }d(\Pi u, \Pi v)\le \frac 12 d(u,v),$$
for any $u,v\in D_\ep$ with $\ep\le \ep_0$. Thus the map $\Pi$ is a
contraction map on the complete set $D_\ep$, and
 the fixed point $u\in D_\ep$ of the map $\Pi$ gives the required
unique almost global solution.

\section{Discussion}
In this section, we consider the generalizations of
Theorem~\ref{thm-almost} to the general problems. To begin, we
consider the following Cauchy problem \beeq\label{SLW-2} \Box u
=\sum_{|\al|=p} P_\al(u) (\pa u)^\al \equiv N(u)\eneq on
$\R\times\R^{n}$, with Cauchy data at time $t=0$ \beeq\label{data-2}
u(0,x)=u_0\in H^s, \ \pt u(0,x) = u_1 \in H^{s-1}\ . \eneq

Before presenting our well posed results, we give a brief history of
the study of these problems. In the case of classical $C_0^\infty$
initial data with size of order $\ep$, the global or almost global
existence can be proved for $p\ge 1+\frac{2}{n-1}$ (global if
$p>1+\frac{2}{n-1}$ and almost global with lifespan $T_\ep \ge
\exp(c\epsilon^{-(p-1)})$ if $p=1+\frac{2}{n-1}$), see Sogge
\cite{So08}. When $p=1+\frac{2}{n-1}$, the lifespan $T_\ep$ is also
sharp for the problem with nonlinearity $|\pt u|^p$ (Zhou
\cite{Zh01}).

As in Theorem \ref{thm-almost}, our object here is to prove the
corresponding results with low regularity. We will consider only the
case $p\ge 3$, which can be easily handled by the linear Strichartz
type estimates.

The scaling consideration tells us the critical Sobolev space for
the equation is $\dot H^{s_c}$ with exponent \beeq\label{sc-2}s_c =
\frac{n+2}{2}-\frac{1}{p-1}\ ,\eneq which is then, heuristically, a
lower bound for the range of permissible $s$ such that the problem
\eqref{SLW}-\eqref{data} is well-posed in $C_t H^s_x\cap C^1_t
H^{s-1}_x$. (See e.g. Theorem 2 in \cite{FaWa08} for the ill posed
result with $s<s_c$ and $N(u)=(\pt u)^p$.)

The local well posedness for the problem with low regularity has
been extensively studied (see Ponce-Sideris \cite{PoSi93}, Tataru
\cite{Ta99} and the authors \cite{FaWa05}). It is proved that the
problem is local well posed in $C_t H^s_x\cap C^1_t H^{s-1}_x$ with
$s>s_c$, if $p\ge 1+ \frac{4}{n-1}$. When $p<1+ \frac{4}{n-1}$,
there is another mechanism due to Lorentz invariance such that the
problem is not well posed in $C_t H^s_x\cap C^1_t H^{s-1}_x$ with
$s=s_c+\ep$ for any $\ep>0$ small enough (see e.g. Lindblad
\cite{Ld93}). In this case ($2\le p<1+ \frac{4}{n-1}$), the local
well posedness is proved for $s>\frac{n+5}4$.

We will prove the corresponding results in the Sobolev spaces $H^s$
with regularity $s\ge s_c$ (for $p> \max(1+ \frac{4}{n-1}, 3)$) or
$s> s_c$ (for $p= \max(1+ \frac{4}{n-1}, 3)$ and $n\neq 3$). For the
remained cases, we will need to use the Sobolev space with certain
angular regularity $b>0$ to establish the global existence. Now we
are ready to state our existence results.

\begin{thm}\label{thm-Glob}
Let $n\ge 3$, $p\ge 3$ and $(n,p)\neq (3,3)$. Then if $p> 3$, the
problem \eqref{SLW-2}-\eqref{data-2} is critically local well-posed
in $C_t H^s_x\cap C^1_t H^{s-1}_x$ such that $\pa u\in L_t^{p-1}
L_x^\infty$, for any $s\ge s_c$. Moreover, if the data is
sufficiently small in $H^s\times H^{s-1}$, then the solution is
global. Furthermore, if $n\ge 4$, $p=3$, we have the same result for
any $s>s_c$.
\end{thm}
\begin{rem}\label{rem}
  The reason we exclude $n=2$ in this Theorem is that we do not have
the Sobolev embedding  $$\dot H^s\cap \dot H^1\subset L^\infty, $$
which is true only for $n\ge 3$. This is also the reason we need to
use Lemma \ref{thm-bound-improved} to give the bound of $u$ in the
proof of Theorem \ref{thm-almost}. However, for the case we have
$P_\al(u)$ are constants $C_\al$, we can still have similar results
for $n=2$. Precisely, the problem is critical local well posed in
$C_t H^s\cap C^1_t H^{s-1}$ for $p>5$ and $s\ge s_c$, and for $p=5$
and $s>s_c$.  It seems interesting to see whether the general
problem \eqref{SLW-2} with $P_\al(u)\neq C_\al$ admits global
solution in this setting.
\end{rem}

The remained case is now $n=2$ or $n=p=3$. In the case $n=p=3$,
global result with small data in $H^{s_c, b}_\theta$ with any $b>0$
has been proved in \cite{MaNaNaOz05}, by proving the angular
Strichartz estimates
$$\|\pa u\|_{L^2_t L^\infty_{|x|}L^r_\theta(\R\times \R^3)}\le C_r
( \|u_0\|_{H^2(\R^3)}+ \|u_1\|_{H^1(\R^3)}),$$
for any $r<\infty$.

Now let us deal with the case $n=2$ and $p\ge 4$. As we have seen in
the Remark \ref{rem}, we are able to prove global well posedness in
$C_t H^s_x\cap C^1_t H^{s-1}_x$ for small data when
$P_\al(u)=C_\al$. However, in view of Lemma
\ref{thm-bound-improved}, it seems difficult to prove the
corresponding result for the general case.

Instead, recall that in the study of the wave equations with
nonlinearities involving $u$, we usually prove global existence by
requiring that the initial data belongs to homogeneous Sobolev
spaces $\dot H^{\gamma}\times \dot H^{\gamma-1}$ with $\gamma\in
(0,1)$ (see e.g. \cite{FaWa10} and Metcalfe-Sogge \cite{MeSo10}).
Inspired by these results, we can remedy this difficulty by
requiring that the second initial data $g\in \dot H^{-\delta}$ in
addition.

Now we can state our results for the case $n=2$.
\begin{thm}\label{thm-Glob-n=2}
Let $n= 2$, $p\ge 5$ and $s\ge s_c$ (with $s>s_c$ when $p=5$).
Consider the equation \eqref{SLW-2}-\eqref{data-2} with fixed
$\delta\in (0,1)$, there exist constants $\ep_0>0$ and $C>0$ such
that we have an unique global solution $u\in C_t H^{s}\cap C_t^1
H^{s-1}$ with
$$\|\pa_{t,x} u\|_{L_t^\infty H^{s-1}\cap L^{p-1}_t L_x^\infty}+\ep \|\pa_{t,x} u\|_{
L_t^\infty \dot H^{-\delta}}\le C\ep,$$ whenever the initial data
satisfies
$$\|\nabla u_0\|_{H^{s-1}}+\|u_1\|_{H^{s-1}}+
\ep \|u_0\|_{\dot H^{1-\delta}}+\ep \|u_1\|_{\dot H^{-\delta}}\le \ep$$
with $\ep\le \ep_0$.
\end{thm}
In the remained case with $p=4$, as in Theorem \ref{thm-almost}, we
need to use Sobolev space with angular regularity as solution space.
\begin{thm}\label{thm-Glob-Angu}
Let $n= 2$, $p\ge 4$, $s\ge s_c$ and $b>1/2$. Consider the equation
\eqref{SLW-2}-\eqref{data-2} with fixed $\delta\in (0,1)$, there
exist constants $\ep_0>0$ and $C>0$ such that we have an unique
global solution $u\in C_t H^{s,b}_\theta\cap C_t^1 H^{s-1,b}_\theta$
with
$$\|\pa_{t,x} u\|_{L_t^\infty H^{s-1,b}_\theta\cap L^{p-1}_t L_{|x|}^\infty H^b_\theta}+\ep \|\pa_{t,x} u\|_{
L_t^\infty \dot H^{-\delta,b}_\theta}\le C\ep,$$ whenever the
initial data satisfies
$$\|\nabla u_0\|_{H^{s-1,b}_\theta}+\|u_1\|_{H^{s-1,b}_\theta}+
\ep \|u_0\|_{\dot H^{1-\delta,b}_\theta}+\ep \|u_1\|_{\dot
H^{-\delta,b}_\theta}\le \ep$$ with $\ep\le \ep_0$.
\end{thm}

\begin{rem}
For the case $p=2$ with $P_\al(u)=C_\al$, there are some related
works of Sterbenz \cite{St04} ($n\ge 6$), \cite{St07} ($n=4$).
\end{rem}

\subsection{Critical Local Well Posedness for $n\ge 3$}
Recall the classical Strichartz estimates (see \cite{FaWa06},
\cite{KeTa98} for example), for the solution to the equation $\Box
u=0$ with initial data $(u_0,u_1)$, we have
\beeq\label{5-Stri}\|\pa_{t,x} u\|_{L^{p-1}_t
L^\infty_{x}}+\|\pa_{t,x} u\|_{L^\infty_t H^{s-1}}\le C (\|\nabla
u_0\|_{H^{s-1}}+\|u_1\|_{H^{s-1}}),\eneq if $p> \max(1+
\frac{4}{n-1}, 3)$ and $s\ge s_c$. Moreover, the same estimates are
true if $p= \max(1+ \frac{4}{n-1}, 3)$, $s>s_c$ with $n\neq 3$.

Let $$d^T_1(u,v)=\|\pa_{t,x} (u-v)\|_{L^{p-1}_t
L^\infty_{x}([0,T]\times \R^n) },$$
$$d_2(u,v)=\|\pa_{t,x} (u-v)\|_{ L^\infty_t H^{s-1}}.$$
We define the solution space to be the complete domain with $\ep\ll
1$ to be chosen later,
\begin{eqnarray*} &  S=&\{u\in C_t H^{s} \cap C_t^1
H^{s-1}: d_1^T(u,0)\le 10 C \ep, \\ &&d_2(u,0)\le 10 C (\|\nabla
u_0\|_{H^{s-1}}+\|u_1\|_{H^{s-1}}):=M\},
\end{eqnarray*}
with metric $d(u,v)=d_1(u,v)+d_2(u,v)$.

For any $u\in S$, we denote $\Pi u$ to be the solution to the linear
equation
$$\Box \Pi u=P_\al(u)(\pa u)^\al$$ with data $(u_0,u_1)$.  Then, using \eqref{5-Stri}, we have
$$  d_1^T(\Pi 0, 0)+ d_2(\Pi 0,0) \le C (\|\nabla u_0\|_{H^{s-1}}+\|u_1\|_{H^{s-1}}),$$
and hence $\lim_{T\rightarrow 0} d_1^T(\Pi 0, 0) =0$. Thus given
$\ep$ small enough, we can choose $T=T(\ep)>0$ small enough (and
$T=\infty$ when the initial data is small enough) such that
$$d_1^T(\Pi 0, 0)\le C \ep .$$ This implies that $\Pi 0\in S$.

Observe that we have $$\|u\|_{L^\infty_{t,x}\cap L^\infty_t \dot
H^{s-1,\frac n{s-1}}}\les \|u\|_{L^\infty_t \dot H^1}+
\|u\|_{L^\infty_t \dot H^{s}}\les d_2(u,0)$$ if $n\ge 3$. If $u\in
S$ and $\ep>0$ is small enough, then by \eqref{5-Stri}, Duhamel's
principle and fractional Leibniz rule, we have
\begin{eqnarray*}
  d_2(\Pi u, 0) &\le& C (\|\nabla u_0\|_{H^{s-1}}+\|u_1\|_{H^{s-1}})+ C
  \|P_\al(u)
  (\pa u)^\al\|_{L^1_t H^{s-1}}\\
&\le & M/10+ \tilde C(M) \|\pa u\|_{L^{p-1}_t L^\infty_{x}}^{p-1}  \|\pa u\|_{L^\infty_t H^{s-1}}\\
&\le & M/10+ \tilde C(M) (10 C \ep)^{p-1} M \le M,
\end{eqnarray*} and similarly
\begin{eqnarray*}
  d_1^T(\Pi u, 0) &\le& d_1^T(\Pi u-\Pi 0, 0)+d_1^T(\Pi 0,0)\\
&\le &   C \|P_\al(u)(\pa u)^\al\|_{L^1_t H^{s-1}} + C \ep\\
&\le & C \ep+ \tilde C(M) \|\pa u\|_{L^{p-1}_t L^\infty_{x}}^{p-1}  \|\pa u\|_{L^\infty_t H^{s-1}}\\
&\le &  C \ep+ \tilde C(M) (10 C \ep)^{p-1} M \le 10 C \ep.
\end{eqnarray*}
This means that $\Pi u\in S$.

Moreover, for any $u,v\in S$, if $\ep\ll 1$ sufficiently small, we
have
\begin{eqnarray*}
  d(\Pi u, \Pi v) &\le& C \|P_\al(u)(\pa u)^\al-P_\al(v)(\pa v)^\al\|_{L^1_T H^{s-1}}\\
&\le& C \|P_\al(u)((\pa u)^\al-(\pa v)^\al)\|_{L^1_T
H^{s-1}}+C\|(P_\al(u)-P_\al(v))(\pa v)^\al\|_{L^1_T
H^{s-1}}\\
&\le & \tilde C(M) (d_1^T( u,0)^{p-2}+ d_1^T( v,0)^{p-2}) (d( u,0)+ d( v,0)) d(u,v)\\
&&+\tilde C(M) (d_1^T( v,0))^{p-1}) d_2( v,0) d_2(u,v)\\
 &\le & 4 M \tilde C(M) (10 C \ep)^{p-2}    d(u,v)+M \tilde C(M) (10 C \ep)^{p-1}    d(u,v) \le \frac{1}{2}
d(u,v).
\end{eqnarray*} Thus we prove that the map $\Pi$ is a contraction
map on the complete set $S$, and then the fixed point $u\in S$ of
the map $\Pi$ gives the required unique solution.

\subsection{Global existence with small data for $n=2$ and $p\ge 4$}
In this subsection, we prove Theorem \ref{thm-Glob-n=2} and Theorem
\ref{thm-Glob-Angu}. The proofs of both Theorems use essentially the
same idea. So for simplicity, we give only the harder proof for
Theorem \ref{thm-Glob-Angu}, which involves the Sobolev space with
angular regularity.

Recall that by Theorem \ref{thm-Strichartz}, for the solution to the
equation $\Box u=0$ with initial data $(u_0,u_1)$, we have
\beeq\label{eq-Stri-end-2}\|\pa_{t,x} u\|_{L^{p-1}_t L^\infty_{|x|}
H^{b}_\theta \cap L^\infty_t H^{s-1,b}_\theta} \le C \|\nabla
u_0\|_{H^{s-1,b}_\theta}+\|u_1\|_{H^{s-1,b}_\theta} := C \ep \ \eneq
for some constant $C>1$. Note also that the initial data satisfy
$$ \|u_0\|_{\dot
H^{1-\delta,b}_\theta}+\|u_1\|_{\dot H^{-\delta,b}_\theta}\le 1.$$

Let $$d_1(u,v)=\|\pa_{t,x} (u-v)\|_{L^{p-1}_t
L^\infty_{|x|}H^{b}_\theta \cap L^\infty_t H^{s-1,b}_\theta}\ ,$$
$$d_2(u,v)=\ep \|\pa_{t,x} (u-v)\|_{ L^\infty_t \dot H^{-\delta,b}_\theta}\ .$$
We define the solution space to be the complete domain with $\ep\ll
1$ to be chosen later, $$S=\{ u\in C_t H^{s,b}_\theta\cap C_t^1
H^{s-1,b}_\theta: d_1(u,0)+d_2(u,0)\le 10 C \ep\},
$$
with metric $d(u,v)=d_1(u,v)+d_2(u,v)$.

As before, for any $u\in S$, we denote $\Pi u$ to be the solution to
the linear equation
$$\Box \Pi u=P_\al(u)(\pa u)^\al$$ with initial data $(u_0,u_1)$.  Then, using \eqref{eq-Stri-end-2}, we have
$$  d_1(\Pi 0, 0)\le C \ep.$$
Moreover, by the energy estimates, we know that
$$  d_2(\Pi 0, 0)\le \ep (\|u_0\|_{\dot
H^{1-\delta,b}_\theta}+\|u_1\|_{\dot H^{-\delta,b}_\theta})\le
C\ep.$$ This implies that $\Pi 0\in S$.

Observe that we have $$\|u\|_{L^\infty_{t,|x|}L^2_\theta \cap
L^\infty_t \dot H^{1}}\les\|u\|_{L^\infty_{t,x}\cap L^\infty_t \dot
H^{1}}\les \|u\|_{L^\infty_t \dot H^{1-\delta}}+ \|u\|_{L^\infty_t
\dot H^{s}}$$ since $s>1=\frac n2$. Then for $u\in S$, we have
\beeq\label{4.6}\|u\|_{L^\infty_{t,|x|}H^b_\theta \cap L^\infty_t
\dot H^{1,b}_\theta}\les 1\ .\eneq

 If $u\in S$ and $\ep>0$
is small enough, then by \eqref{eq-Stri-end-2}, Duhamel's principle
and fractional Leibniz rule (Lemma \ref{thm-Leibniz}), we have
\begin{eqnarray*}
  d_1(\Pi u, 0) &\le& C \ep + C
  \|P_\al(u)
  (\pa u)^\al\|_{L^1_t H^{s-1,b}_\theta}\\
&\le & C \ep + C_1 \|\pa u\|_{L^{p-1}_t L^\infty_{|x|} H^b_\theta}^{p-1}
\|\pa u\|_{L^\infty_t H^{s-1,b}_\theta}\\
&\le & C\ep + C_1 d_1(u,0)^p \le 2 C\ep\ .
\end{eqnarray*}

To give the estimate of $d_2$, we need to use the following
generalized Strichartz estimates (see Remark \ref{rem-3} or
Proposition 1.2 in \cite{SmSoWa09})
\beeq\label{4.7}\|e^{-itP}f\|_{L^q_t L^r_{|x|}L^2_\theta}\le C_{q,r}
\|f\|_{\dot H^{1-\frac 1q-\frac 2r}}\eneq if $\frac 1q+\frac
1r<\frac 12$. We will only use the special case when
$$q=\frac{p-2}\delta\textrm{ and }1-\frac 1q - \frac 2r=\delta\ .$$
This choice of $(q,r)$ satisfies the relation
$\frac{p-1}{q'}=p-\frac 2{r'}$. Note here that we have the
admissible condition for \eqref{4.7} with this choice of $(q,r)$ for
small enough $\delta>0$ if and only if $p>3$.

Recall that by Sobolev embedding, the following embedding estimates
are true
$$L^{2/(1+\delta)}_{|x|}L^2_\theta\subset L^{2/(1+\delta)}\subset \dot H^{-\delta}
\textrm{ and } \dot H^{1-\delta}\subset L^{2/\delta}\ .$$ Thus by
\eqref{4.6}, \eqref{4.7}, duality, H\"{o}lder inequality and the
fact that $H^{b}_\theta$ is an algebra under multiplication when
$b>1/2$, we have
\begin{eqnarray*}
  d_2(\Pi u, 0) &\le& d_2(\Pi u-\Pi 0, 0)+d_2(\Pi 0,0)\\
&\le &   \ep \|(P_\al(u)-P_\al(0))(\pa u)^\al\|_{L^1_t \dot
H^{-\delta,b}_\theta}+ \ep \|P_\al(0)(\pa u)^\al\|_{ L^{q'}_t
L^{r'}_{|x|}H^b_\theta} + C \ep\\
&\le & C \ep+\ep \|(P_\al(u)-P_\al(0))(\pa u)^\al\|_{ L^1_t
L^{2/(1+\delta)}_{|x|}H^b_\theta}+ \ep  |P_\al(0)| \|(\pa u)^\al\|_{
L^{q'}_t L^{r'}_{|x|}H^b_\theta}\\
&\le & C \ep+\ep \|P_\al(u)-P_\al(0)\|_{L^\infty_t
L^{2/\delta}_{|x|}H^b_\theta}\| (\pa u)^\al\|_{
L^1_t L^{2}_{|x|}H^b_\theta}\\
&&+\ep  |P_\al(0)| \|\pa u\|_{L^\infty_t L^2_{|x|}H^b_\theta}^{2/r'}
\|\pa u\|_{L^{p-1}_t L^\infty_{|x|}H^b_\theta}^{p-2/r'}\\
&\le & C \ep+C_2 (d_2(u,0)+\ep)  d_1(u,0)^p\\
&\le & 2 C \ep\ .
\end{eqnarray*}
This means that $\Pi u\in S$.

Moreover, for any $u,v\in S$, if $\ep\ll 1$ sufficiently small, we
have
\begin{eqnarray*}
  d_1(\Pi u, \Pi v) &\le& C \|P_\al(u)(\pa u)^\al-P_\al(v)(\pa v)^\al\|_{L^1_T H^{s-1,b}_\theta}\\
&\le& C \|P_\al(u)((\pa u)^\al-(\pa v)^\al)\|_{L^1_T
H^{s-1}}+C\|(P_\al(u)-P_\al(v))(\pa v)^\al\|_{L^1_T
H^{s-1}}\\
&\le & C_3 (d_1( u,0)^{p-1}+ d_1( v,0)^{p-1}) d(u,v)\\
&&+C_3 (d_1( v,0))^{p}
\|P_\al(u)-P_\al(v)\|_{L^\infty_{t,|x|}H^b_\theta \cap L^\infty_t
\dot H^{1,b}_\theta}\\
&\le & C_3 (d_1( u,0)^{p-1}+ d_1( v,0)^{p-1}) d(u,v) + \tilde C_3
\ep^{-1} (d_1( v,0))^{p} d(u,v)\\
 &\le & (2^{p-1} C_3 (10 C \ep)^{p-1}+ \tilde C_3 \ep^{-1} (10 C \ep)^{p}) d(u,v)
 \le \frac{1}{2}
d(u,v).
\end{eqnarray*}
As in the proof of $d_2(\Pi u,0)$, we have
\begin{eqnarray*}
  d_2(\Pi u, \Pi v) &\le&
\ep \|(P_\al(u)-P_\al(0))((\pa u)^\al-(\pa v)^\al)\|_{L^1_t
  L^{2/(1+\delta)}_{|x|}H^b_\theta}\\
  &&+\ep\|(P_\al(u)-P_\al(v))(\pa v)^\al\|_{L^1_t
  L^{2/(1+\delta)}_{|x|}H^b_\theta}\\
  &&+\ep \|P_\al(0)((\pa u)^\al-(\pa v)^\al)\|_{L^{q'}_t
L^{r'}_{|x|}H^b_\theta}
  \\
&\le & C_4 d_2(u,0) (d_1( u,0)^{p-1}+ d_1( v,0)^{p-1}) d_1(u,v)\\
&&+C_4 d_2(u,v) d_1( v,0)^{p}\\
&&+C_4 \ep ((d_1( u,0))^{p-1}+(d_1( v,0))^{p-1}) d_1(u,v)\\
 &\le & 10 C_4 (10 C \ep)^{p}    d(u,v)\le \frac{1}{2}
d(u,v).
\end{eqnarray*}

 Thus we prove that the map $\Pi$ is a contraction
map on the complete set $S$, and then the fixed point $u\in S$ of
the map $\Pi$ gives the required unique solution.

\medskip
\noindent{\sl Acknowledgement.} The second author would like to
thank Professor Christopher Sogge for helpful discussion on the
Strichartz estimates. The authors would also like to thank Ting
Zhang for useful communication on the fractional Leibniz rule.


\end{document}